\documentclass[12pt]{article}%
\usepackage{amsfonts}
\usepackage{amssymb}
\usepackage{amsmath}
\usepackage{graphicx}%
\providecommand{\U}[1]{\protect\rule{.1in}{.1in}}
\providecommand{\U}[1]{\protect\rule{.1in}{.1in}}

\begin{document}

\title{Schur's exponent conjecture --- counterexamples of exponent $5$ and exponent
$9$}
\author{Michael Vaughan-Lee}
\date{August 2020}
\maketitle

\section{Introduction}

There is a long-standing conjecture attributed to I. Schur that if $G$ is a
finite group with Schur multiplier $M(G)$ then the exponent of $M(G)$ divides
the exponent of $G$. It is easy to show that this is true for groups $G$ of
exponent 2 or exponent 3, but it has been known since 1974 that the conjecture
fails for exponent 4. Bayes, Kautsky and Wamsley \cite{bayes74} give an
example of a group $G$ of order $2^{68}$ with exponent 4, where $M(G) $ has
exponent 8. (Bayes, Kautsky and Wamsley are heros of the early days of
computing with finite $p$-groups.) However the truth or otherwise of this
conjecture has remained open up till now for groups of odd exponent, and in
particular it has remained open for groups of exponent 5 and exponent 9. For a
survey article on Schur's conjecture see Thomas \cite{viji}.

In this note I give an example of a four generator group $G$ of order
$5^{4122}$ with exponent 5, where the Schur multiplier $M(G)$ has exponent 25,
and an example of a four generator group $A$ of order $3^{11983}$ and exponent
9, where the Schur multiplier $M(A)$ has exponent 27. Very likely the reason
that similar examples have not been found up till now is that computing the
Schur multipliers of groups of this size is right on the edge of what is
possible with today's computers.

We define the group $G$ as follows. First we let $H$ be the four generator
group with presentation
\[
\langle a,b,c,d\,|\,[b,a]=[d,c]\rangle,
\]
and then we let $G$ be the largest quotient of $H$ with exponent 5 and
nilpotency class 9. Let $F$ be the free group of rank 4, with free generators
$a,b,c,d$, and let $M$ be the normal closure in $F$ of 
$\{g^{5}\,|\,g\in F\}\cup\{[b,a][c,d]\}$. Then $G=F/R_G$ where $R_G=M\gamma_{10}(F)$. 
The group $F/[R_G,F]$ is a central extension of $G$, and the Schur multiplier $M(G)$ 
is $(R_G\cap F^{\prime })/[R_G,F]$. Clearly $[b,a][c,d]\in R_G\cap F^{\prime}$, and 
we show that $G$ is a counterexample to the Schur exponent conjecture by showing 
that $([b,a][c,d])^{5}\notin\lbrack R_G,F]$.

The group $A$ is defined similarly. It is the largest quotient of%
\[
\langle a,b,c,d\,|\,a^{3},\,b^{3},\,c^{3}\,,d^{3},\,[b,a]=[d,c]\rangle
\]
with exponent 9 and nilpotency class 9. So if we let $N$ be the normal closure in $F$
of
\[
\{g^{9}\,|\,g\in F\}\cup\{a^3,b^3,c^3,d^3,[b,a][c,d]\}
\]
then $A=F/R_A$ where $R_A=N\gamma_{10}(F)$, and the Schur multiplier $M(A)$ is
$(R_A\cap F^{\prime })/[R_A,F]$. We show that $A$ is a counterexample to the Schur
exponent conjecture by showing that $([b,a][c,d])^{9}\notin\lbrack R_A,F]$.

I was led towards these examples after a fruitful correspondence with Viji
Thomas. He wrote to me saying that he was investigating the groups $R(d,5)$
for various $d$. (Here $R(d,5)$ is the largest finite quotient of the $d$
generator Burnside group of exponent 5, $B(d,5)$.) He mentioned that the Schur
exponent conjecture was still open for groups of exponent 5, but that he could
prove that the Schur multipliers of $R(2,5)$ and $R(3,5)$ have exponent 5. He
wondered if I knew what the nilpotency class of $R(4,5)$ is. It is known that
the class of $R(d,5)$ is at most $6d$ (see \cite{havasnvl90}), so that the
class of $R(4,5)$ is at most 24. He said that if in fact the class is less
than 24 then he might be able to prove that the Schur multiplier of $R(4,5)$
has exponent 5. It seems quite likely that the class of $R(4,5)$ is less than
24 since the class of $R(3,5)$ is 17. (The class of $R(2,5)$ is 12.) But I was
unable to help him on this point since as far as I know the class of $R(4,5)$
remains undetermined. Out of interest I computed the Schur multiplier of
$R(2,5)$ --- it is elementary abelian of order $5^{31}$. Detailed information
from this computation led me to conclude that any exponent 5 counterexample to
Schur's exponent conjecture would need to have class at least 9 and would need
at least 4 generators. This detailed information also showed that the Schur
multiplier of $R(4,5)/\gamma_{10}(R(4,5))$ has exponent 5. So if we want to 
find a class 9 quotient $G$ of $R(4,5)$ with Schur multiplier with exponent 
greater than 5, then $G$ needs to satisfy a relation $r=1$ where $r$ is a 
product of commutators of weight at least 2, and where $r=1$ is not a consequence 
of fifth power relations. This is what led me to consider the relation $[b,a]=[d,c]$.

In the next section I show one way of computing the Schur multiplier of $R(2,5)$, 
and then in Section 3 I show how to compute the central extension $F/[R_G,F]$ of $G$. 
Success with the group $G$ led me to investigate the group $A$ of exponent 9, and this
calculation is described in Section 4. Finally, in Section 5 I speculate on
possible counterexamples in other exponents.

\section{The Schur multiplier of $R(2,5)$}

The group $R(2,5)$ has order $5^{34}$ and nilpotency class 12. You can verify
this in \textsc{Magma} \cite{boscan95} by entering
\[
\mathrm{P:=pQuotient(FreeGroup(2),5,0:Exponent:=5,Print:=1);}
\]
The $p$-covering group of the class 11 quotient of $R(2,5)$ has order $5^{65}$,
and so as a class 12 group $R(2,5)$ has a presentation with 31 fifth powers as
relators. If we take generators $a,b$ for $R(2,5)$ then a suitable set of
relators is $\{u^{5}\,|\,u\in U\}$ where $U$ consists of the elements%
\begin{align*}
&  a,b,ab,a^{2}b,ab^{2},a^{3}b,a^{2}b^{2},ab^{3},a^{4}b,a^{3}b^{2}%
,a^{2}bab,a^{3}bab,a^{2}bab^{2},a^{2}b^{2}ab,abab^{3},\\
&  a^{4}bab,a^{3}bab^{2},a^{2}bab^{3},a^{2}b^{2}ab^{2},a^{4}bab^{2}%
,a^{3}ba^{2}b^{2},a^{3}bab^{3},a^{2}babab^{2},a^{2}bab^{4},\\
&  a^{2}b^{2}ab^{3},a^{4}ba^{2}b^{2},a^{4}bab^{3},a^{3}ba^{2}bab,a^{3}%
babab^{2},a^{3}bab^{4},a^{2}ba^{2}bab^{2}.
\end{align*}
So if we let $F_2$ be the free group of rank 2 generated by $a,b$ and let $K$ be
the normal closure in $F_2$ of $\{u^5\,|\,u\in U\}$ then $R(2,5)=F_2/R$ where
$R=K\gamma_{13}(F_2)$. (You can use the $p$Quotient algorithm in \textsc{Magma}
to verify that $F_2/R$ has order $5^{34}$.) Let $S$ be the central extension
$F_2/[R,F]$ of $R(2,5)$. Then $S$ is the class 13 quotient of the group
with presentation
\[
\langle a,b\,|\,\{[u^{5},v]\,|\,u\in U,\,v\in\{a,b\}\}\rangle.
\]
In \textsc{Magma} you can compute a PC-presentation for $S$
using the nilpotent quotient algorithm. If we let $T$ be the subgroup $\langle
u^{5}\,|\,u\in U\rangle\gamma_{13}(S)$ of $S$, then $S/T$ is isomorphic to
$R(2,5)$ and the Schur multiplier of $R(2,5)$ is $T\cap S^{\prime}$. 
As mentioned above, the Schur multiplier is elementary abelian of order $5^{31}$.

One important observation is that $[b,a]^{5}\in\gamma_{10}(S)$. This implies
that if $P$ is a finite group of exponent 5 with class less than 9 then the
derived group of any central extension of $P$ has exponent 5. And this implies that
the Schur multiplier $M(P)$ has exponent 5. In fact detailed examination shows
that in $S$ we can express $[b,a]^{5}$ as a product of commutators
$[x_{1},x_{2},\ldots,x_{k}]\,$\ $(k\geq10)$ where $x_{1},x_{2},\ldots,x_{k}%
\in\{a,b\}$ and where $a$ and $b$ both occur at least 5 times in the sequence
$x_{1},x_{2},\ldots,x_{k}$. This implies that if $H$ is a central extension of any
group of exponent 5, and if $c\in H$ is a commutator of weight $k>1$ then
$c^{5}\in\gamma_{5k}(H)$.

\section{Computing the group $F/[R_G,F]$}

As stated in the Introduction, we let $G$ be the largest exponent 5, class 9,
quotient of $\langle a,b,c,d\,|\,[b,a]=[d,c]\rangle$. We write $G=F/R_G$ where
$F$ is the free group of rank 4 with free generators $a,b,c,d$. We want to 
compute $F/[R_G,F]$.

We can use the $p$Quotient algorithm in \textsc{Magma} with parameter
\textquotedblleft Exponent:=5\textquotedblright\ to compute a PC-presentation
for $G$. This takes about 4 minutes of CPU-time. (All timings are for programs
run in \textsc{Magma} V2.19-10 running on a desktop computer with 16GB of RAM
and an Intel Core i7-4770CPU@3.40GHz$\times$8 processor. This is quite an old
version of \textsc{Magma} and John Cannon keeps telling me that I really ought
to upgrade to the latest version.) The calculation shows that $G$ has order
$5^{4122}$. As in our computation of the Schur multiplier of $R(2,5)$ we need to
find a finite set of fifth powers which together with the relation
$[b,a]=[d,c]$ define $G$ as a class 9 group. There is a theorem of Higman
\cite{higman59} which implies that if $G$ is nilpotent of class $c$ then $G$
has exponent dividing $n$ provided $g^{n}=1$ for all words of length at most
$c$ in the generators of $G$. So I generated a list of all words of length at
most 9 in the generators $a,b,c,d$ of $G$. There are some obvious redundancies
in this list. For example $ab$ is conjugate to $ba$ so that the relation
$(ab)^5=1$ is equivalent to the relation $(ba)^5=1$, and $ba$ is redundant. More
generally, if $x_{1},x_{2},\ldots,x_{k}\in\{a,b,c,d\}$ then $x_{1}x_{2}\ldots
x_{k}$ is conjugate to $x_{2}\ldots x_{k}x_{1}$ and so we can discard any word
which is lexicographically greater than any of its cyclic conjugates. We can
also discard any word which contains a subword $a^{5}$, $b^{5}$, $c^{5}$ or
$d^{5} $. This left me with a list of 39564 words in the generators $a,b,c,d$.
In principle you could use this list to compute $F/[R_G,F]$, but
the computation would probably take a month or more of CPU-time. The
$p$-covering group of the class 8 quotient of $G$ has order $5^{7044}$ so we
need 2921 fifth powers (together with the relation $[b,a]=[d,c]$) to define
$G$ as a class 9 group. I reduced my long list of fifth powers to a list of
2921 fifth powers%
\[
a^{5},b^{5},c^{5},d^{5},(ab)^{5},\ldots,(ab^{4}cbcd)^{5}%
\]
as follows. First I computed the class 5 quotient $K$ of%
\[
\langle a,b,c,d\,|\,[b,a]=[d,c],a^{5},b^{5},c^{5},d^{5}\rangle.
\]
This group $K$ has order $5^{214}$, whereas the class 5 quotient of $G$ has
order $5^{162}$. So $|K^{5}|=5^{52}$. I systematically built up the subgroup
$K^{5}$, starting with the trivial subgroup $L=\{1\}$ and adding in one fifth
power at a time to $L$ from my long list of fifth powers, till $L$ had order
$5^{52}$. By keeping track of which fifth powers increased the order of $L$, I
was able to obtain a list of 52 fifth powers which together with the relations
$[b,a]=[d,c]$, $a^{5}=1$, $b^{5}=1$, $c^{5}=1$, $d^{5}=1$ define the class 5
quotient of $G$. Next I computed the class 6 quotient $M$ of the group
satisfying these 52 fifth power relations in addition to the relations
$[b,a]=[d,c]$, $a^{5}=1$, $b^{5}=1$, $c^{5}=1$, $d^{5}=1$. Then I found a
minimal set of fifth powers from the long list of fifth powers which generate
$M^{5}$. And so on, up to class 9. Tedious, but straightforward enough.

I now had a list $U$ of 2921 words in $a,b,c,d$ with the property that $G$ is
the class 9 quotient of the group with generators $a,b,c,d$ and relations
\[
\{u^{5}=1\,|\,u\in U\}\cup\{[b,a]=[d,c]\}.
\]
As a check I ran the $p$Quotient algorithm up to class 9 on these generators
and relations. This took 12 minutes of CPU-time. (I can send the list $U$ to
any reader who is interested in following this up.)

Now let $V=\{u^{5}\,|\,u\in U\}\cup\{[b,a][c,d]\}$, and let%
\[
W=\{[v,w]\,|\,v\in V,\,w\in\{a,b,c,d\}\}.
\]
Then $F/[R_G,F]$ is the class 10 quotient of $\langle a,b,c,d\,|\,W\rangle $.
Call this quotient $S$.

The natural approach would be to use the nilpotent quotient algorithm
to compute a PC-presentation for $S$, as I did when computing the Schur multiplier
of $R(2,5)$. But a computation with the nilpotent quotient algorithm would
have taken months of CPU-time (even if it ever completed). I tried using the
nilpotent quotient algorithm to compute $G$ using the presentation with 2921
fifth powers, and I had to kill the job when it had still not completed after
24 hours. So I used the $p$Quotient algorithm to compute the $p$-class 10
quotient $P$ of $\langle a,b,c,d\,|\,W\rangle$. This took 46 hours of
CPU-time, and showed that $P$ has order $5^{13330}$. Clearly $P$ is a
homomorphic image of $S$ (since $a^{5^{10}}=b^{5^{10}%
}=c^{5^{10}}=d^{5^{10}}=1$ in $P$), but $[b,a][c,d]$ has order 25 in $P$, and
so order at least 25 in $S$. So the Schur multiplier of $G$ has exponent at
least 25. On the other hand we know from the computation of the Schur multiplier of
$R(2,5)$ that $S^{\prime5}\leq\gamma_{10}(S)$, and that $\gamma_{10}(S)$ has
exponent 5. So the Schur multiplier of $G$ has exponent 25.

\section{Computing a quotient of $F/[R_A,F]$}

As stated in the Introduction, we let $A$ be the largest exponent 9, class 9,
quotient of
\[
\langle a,b,c,d\,|\,a^3,b^3,c^3,d^3,[b,a]=[d,c]\rangle .
\]
We write $A=F/R_A$ where
$F$ is the free group of rank 4 with free generators $a,b,c,d$. We want to 
compute $F/[R_A,F]$ (or a suitable quotient of this group).

Let $B$ be the group generated by $a,b,c,d$ with relations%
\[
\{a^{3}=1,\,b^{3}=1,\,c^{3}=1,\,d^{3}=1,\,[b,a][d,c]=1\}\cup\{u^{9}=1\,|\,u\in U\}
\]
where $U$ is the set%
\begin{align*}
& \{ab,ac,ad,bc,bd,cd,a^{2}b,a^{2}c,a^{2}d,abc,abd,acd,b^{2}c,b^{2}d,bcd,\\
& c^{2}d,a^{2}bc,a^{2}bd,a^{2}cd,ab^{2}c,ab^{2}d,abc^{2},abcd,abd^{2}%
,ac^{2}d,acd^{2},b^{2}cd,\\
& bc^{2}d,bcd^{2},a^{2}bcd,ab^{2}cd,abc^{2}d,abcd^{2},a^{2}b^{2}cd,a^{2}%
bc^{2}d,a^{2}bcd^{2}\}.
\end{align*}
Then the class 9 quotient of $B$ has exponent 9, and so is isomorphic to $A$. 
(You can check this in \textsc{Magma} by using the $p$Quotient algorithm to 
compute the class 9 quotient of $B$, and then running the $p$Quotient algorithm 
up to class 9 again, with the extra parameter \textquotedblleft Exponent:=9".) 
So if we let $L$ be the normal closure of
\[
\{u^{9}\,|\,u\in U\}\cup\{a^3,b^3,c^3,d^3,[b,a][c,d]\}
\]
in the free group $F$, then $A=F/R_A$ where $R_A=L\gamma_{10}(F)$

Now let $S$ be the class 10 quotient of the group with generators $a,b,c,d$
and relators%
\[
\{a^{3},b^{3},c^{3},d^{3}\}\cup\left\{  \lbrack x,y]\,|\,x\in \{u^9\,|\,u\in
U\}\cup\{[b,a][c,d]\},\;y\in\{a,b,c,d\}\right\}  .
\]
Then $S$ is a central extension of $A$, and is a proper quotient of the
group $F/[R_A,F]$. It takes the $p$Quotient algorithm in
\textsc{Magma} two minutes to compute $S$, which has order $3^{37170}$.
(Presumably this calculation is so quick compared with the calculation of the
$p$-class 10 quotient of $F/[R_G,F]$ because $A$ has many fewer relations than $G$.)
Unfortunately \textsc{Magma} crashes immediately after completing the
calculation of $S$. It seems to me likely that \textsc{Magma} has a problem
converting $p$Quotient's internal representation of $S$ into a standard \textsc{Magma}
PC-presentation. However \textsc{Magma}'s C version of $p$Quotient is based on
George Havas's original Fortran version \cite{havas2}, and so I recomputed $S$
using George's Fortran code. The computation showed that $[b,a][c,d]$ has
order 27 in $S$, and so order at least 27 in $F/[R_A,F]$.

I notified Eamonn O'Brien, who wrote the $p$Quotient program in \textsc{Magma},
about my problem with \textsc{Magma} crashing. He confirmed that there is
a bug in \textsc{Magma}, even in the latest version. However, with Eamonn's special 
knowledge of his program he was able to use $p$QuotientProcess to confirm my 
Fortran calculation.

So the Schur multiplier $M(A)$ has exponent at least 27. However if we let $T$
be any central extension of a class 9 group of exponent 9 then it is easy to see
that the derived group $T'$ has exponent dividing 27. We proceed as follows.
We can use the nilpotent quotient algorithm to compute the class 10 quotient of
\[
\langle a,b\,|\,\{[u^9,v]\,|\,u\in \{a,b,ab,a^2b,ab^2\},\,v\in \{a,b\}\}\rangle.
\]
The commutator $[b,a]$ has order 27 in this quotient, and so any commutator
in $T$ has order dividing 27. So $T'$ is generated by elements of order at most
27, and has class at most 5. We can use the nilpotent quotient algorithm
to compute the class 5 quotient of
\[
\langle a,b\,|\,\{[u^9,v]\,|\,u\in \{a,b,ab,a^2b,ab^2\},\,v\in \{a,b\}\}\cup 
\{a^{27},b^{27}\}\rangle ,
\]
and $ab$ has order 27 in the quotient. So the product of elements in $T'$ with
order dividing 27 also has order dividing 27. So $T'$ has exponent dividing 27.

All this shows that the Schur multiplier $M(A)$ has exponent 27.

\section{Other exponents?}

It seems certain that there are similar examples for all prime powers greater
than 3. George Havas conjectures that for every prime $p>3$ the largest
exponent $p$, class $2p-1$ quotient of%
\[
\langle a,b,c,d\,|\,[b,a]=[d,c]\rangle
\]
is a counterexample. Certainly it is easy to show that any exponent 7
counterexample must have class at least 13. The problem with computing this
group, even for $p=7$, is not so much that computers nowadays do not have
enough memory or that the calculation would take too long. The problem is
rather that the data structures built into current implementations\ of the $p
$Quotient algorithm never anticipated handling groups of this size. For
example, in George's Fortran program a \textquotedblleft generator exponent
pair\textquotedblright\ $a_{i}^{j}$ is stored as a single 32 bit integer
$2^{16}j+i$, so some adjustment is needed to the data structure if the program
is to be able to handle more that 65535 PC-generators.

\end{document}